\documentclass{amsart}
\usepackage{amssymb,latexsym,enumerate,amscd}
\newtheorem{Thm}[equation]{Theorem}
\newtheorem{Prop}[equation]{Proposition}
\newtheorem{Lem}[equation]{Lemma}
\newtheorem{Cor}[equation]{Corollary}

\theoremstyle{remark}

\newtheorem*{Rem*}{Remark}
\theoremstyle{definition}

\newtheorem*{Not*}{Notation}
\newtheorem*{Def*}{Definition}
\numberwithin{equation}{section}
\DeclareMathOperator{\im}{im}

\DeclareMathOperator{\gr}{gr}
\DeclareMathOperator{\CS}{CS}

\DeclareMathOperator{\coker}{coker}

\DeclareMathOperator{\Mat}{M}

\begin{document}

\date{%
Wed Mar 15 19:01:08 EST 2006}

\title[The Atiyah Conjecture]
{Congruence subgroups and the Atiyah conjecture}

\author[D. R. Farkas]{Daniel R. Farkas}
\address{Department of Mathematics \\
Virginia Tech \\
Blacksburg \\
VA 24061-0123 \\
USA}
\email{farkas@math.vt.edu}
\urladdr{http://www.math.vt.edu/people/farkas/}

\author[P. A. Linnell]{Peter A. Linnell}
\address{Department of Mathematics \\
Virginia Tech \\
Blacksburg \\
VA 24061-0123 \\
USA}
\email{linnell@math.vt.edu}
\urladdr{http://www.math.vt.edu/people/linnell/}

\begin{abstract}
Let $\overline{\mathbb{Q}}$ denote the algebraic
closure of $\mathbb{Q}$ in $\mathbb{C}$.  Suppose $G$ is a
torsion-free group which contains a
congruence subgroup as a normal subgroup of finite index and denote by
$\mathcal{U}(G)$ the $\mathbb{C}$-algebra of
closed densely defined unbounded operators affiliated to the group
von Neumann algebra of $G$.  We prove that there
exists a division ring $D(G)$ such that $\overline{\mathbb{Q}}[G]
\subseteq D(G)
\subseteq \mathcal{U}(G)$.  This establishes some versions of the
Atiyah conjecture for the group $G$.
\end{abstract}



\maketitle

\section{Introduction} \label{Sintroduction}

Let $p$ be a prime and let $d,i$ be positive integers.
Throughout this paper, $\mathbb{Z}_p$ denotes the $p$-adic
integers, $\Mat_d(R)$ denotes the $d \times d$ matrices over
a ring $R$, and $\CS(i,d,p)$ denotes the congruence subgroup
\[
\{A \in \Mat_d(\mathbb{Z}_p) \mid A \equiv I_d \mod p^i\}.
\]
(Here $I_d$ is the identity $d\times d$
matrix.) We will \emph{always} assume that $i \ge 2$ whenever the
prime $p$ is $2$.

We use Passman's notation $R[G]$ for the group ring of $G$ with
coefficients in the ring $R$.
Let $\mathcal{U}(G)$ denote the unbounded linear operators
affiliated to the group von Neumann algebra of $G$.  A brief
description of this object is found at the end of this
introduction.  The main result of this paper is
\begin{Thm} \label{Tmain}
Let $d,u$ be positive integers, let $p$ be a prime,
and let $G$ be a torsion-free group which contains a normal
subgroup of finite index isomorphic to $\CS(u,d,p)$.  Then
there is a division ring $D(G)$ such that $\overline{\mathbb{Q}}[G]
\subseteq D(G) \subseteq \mathcal{U}(G)$.
\end{Thm}

Some years ago, we wrote in \cite{FarkasLinnell89} an exposition of
Moody's and Lazard's contributions to the zero divisor problem.
While the first result has been discussed and exploited extensively
in the literature, Lazard's observation
in \cite{Lazard65}, that the group ring of a
congruence subgroup over the integers is a domain, has not had
nearly as much traction.  Despite our best efforts, this theorem
still has a reputation of being intricate and difficult.  Lazard's
work focuses on an associated graded
ring constructed from finite slices of a congruence subgroup.
In this note, we highlight
the idea that one can detect whether an element of the group ring is a
zero divisor by studying properties of its images in the group rings
of various finite quotient groups.  A direct approach is possible due
to the very clear discussion of filtrations in
characteristic $p$ from \cite[Chapter 3]{Passman77}. (Don Passman's
magnum opus was affectionately referred to as \emph{The Group
Algebraist's Guide to the Galaxy} at the conference celebrating his
achievements.)

For the remainder of the paper, fix positive integers $u$ and $d$
as well as a prime $p$.  Set $\Gamma = \CS(u,d,p)$.
Again, we assume that $u \ge 2$ whenever $p=2$.
(This technical point ensures, among other things, that $\Gamma$ is
torsion-free.) We begin with a purely expository presentation of
\begin{Thm} \label{Tcharp}
If $k$ is a field of characteristic $p$ then $k[\Gamma]$ is a domain.
\end{Thm}
It does not appear that the theorem in this generality can be obtained
immediately from \cite{Lazard65}.  Perhaps it is folklore.
In any event we can now pass to the case of arbitrary fields of
characteristic zero by the following well-known observation:
if $p$ is a fixed prime and $G$ is a group with the property that
$k[G]$ is a domain
for all fields $k$ of characteristic $p$, then $K[G]$ is a domain
for all fields $K$ of characteristic 0.  One can obtain this from
\cite[p.~617]{Passman77} and \cite[Exercises 20 and 21]{Passman77}.
Thus we have
\begin{Thm} \label{Tchar0}
If $K$ is a field of characteristic $0$ then $K[\Gamma]$ is a domain.
\end{Thm}

We shall use the techniques employed in Theorem \ref{Tcharp} to prove
a version of the Atiyah conjecture \cite[\S 10]{Lueck02}.
Let $G$ be an arbitrary group and write
\[
L^2(G) = \{\sum_{g \in G}a_gg \mid a_g \in \mathbb{C}
\text{ and } \sum_{g \in G} |a_g|^2 < \infty \}.
\]
The usual convolution product on $\mathbb {C}[G]$ extends to a module
action $\mathbb {C}[G] \times L^2(G) \to L^2(G)$.  Specifically, if
$\alpha = \sum_{g\in G} a_gg$ and $\beta = \sum_{g\in G} b_gg$,
where $a_g = 0$ for all but finitely many $g$, then
\[
\alpha\beta = \sum_{g \in G}(\sum_{h \in G} a_{gh^{-1}}b_h)g.
\]
Let $\mathcal{B}(L^2(G))$ denote the bounded linear operators on
$L^2(G)$.  We may view elements of $\mathbb{C}[G]$ as bounded
linear operators acting
by left multiplication on $L^2(G)$.  The group von Neumann algebra
$\mathcal{N}(G)$ is the weak closure of $\mathbb{C}[G]$ in
$\mathcal{B}(L^2(G))$.
Let $\mathcal{U}(G)$ denote the ring of unbounded operators
affiliated to the group von Neumann algebra \cite[\S 8.1]{Lueck02}.
Then we have the following inclusions \cite[(8.1), p.~224]{Linnell98}:
\[
\mathbb{C}[G] \subseteq \mathcal{N}(G) \subseteq L^2(G) \subseteq
\mathcal{U}(G).
\]
The relevant operations are, of course, compatible.
One version of the Atiyah conjecture states that if $G$
is a torsion-free group, then there exists a division ring $D$ such
that $\mathbb{C}[G] \subseteq D \subseteq \mathcal{U}(G)$.  Although
we cannot prove this for virtual congruence subgroups, we do
obtain a rational analogue.

\section{Background}

As might be expected in a tribute to Don Passman, we use nontrivial
techniques from both group theory and ring theory.  We also require
heavy doses of topology and analysis.

The exposition begins with rings.  Given a ring $R$, we say a sequence
of ideals $R=I_{0}, I_{1}, I_{2}, \ldots $ is \emph{descending}
provided that
\[
I_{0}\supseteq I_{1}\supseteq I_{2} \supseteq \cdots
\]
and $\bigcap _{t}I_{t}=0$.  It is a \emph{descending filtration} if,
in addition, it satisfies $I_{s}I_{t} \subseteq I_{s+t}$ for all $s$
and $t$.  In this situation, the associated graded ring is
\[
\gr(R)=\bigoplus _{j=0}^{\infty}I_{j}/I_{j+1}.
\]
If $r \in I_{s}\setminus I_{s+1}$ then its image in $\gr(R)$ is the
\emph{bottom symbol} of $r$ and we refer to $s$ as the degree of $r$.
For any group $G$ and any field $k$, the augmentation ideal in the
group ring $k[G]$ is denoted $\omega (G)$.  We adopt the convention
that $\omega (G)^{0}=k[G]$.  If $k$ is a field of prime
characteristic $p$ and $G$ is a finite $p$-group then the powers of
$\omega (G)$ constitute a descending filtration of $k[G]$.  This
entire paper can be regarded as an obsession with this filtration.

Next, we review some notation for groups.
For subgroups $A$ and $B$, we write $[A,B]=
\langle aba^{-1}b^{-1} \mid a \in A,\, b \in B\rangle$,
the commutator group of $A$ and $B$.  For a non-negative integer
$r$, we let $A^r = \langle a^r \mid a \in A\rangle$.

We remind the reader of the fundamental properties of congruence
subgroups.  Clearly
\[
\CS(j,d,p)/\CS(j+1,d,p) \simeq \Mat_{d}(\mathbb{Z}/p\mathbb{Z}),
\]
the group of matrices under addition.  As a consequence, the factor
group is an elementary abelian $p$-group of rank $d^{2}$ for all $j$.
\emph{For the remainder of the paper we set $e=d^{2}$.} An easy
application of the binomial theorem shows that the $p$-power map
induces an isomorphism
\[
\CS(j,d,p)/\CS(j+1,d,p) \rightarrow \CS(j+1,d,p)/\CS(j+2,d,p).
\]

By using a variant of Hensel's Lemma, one proves that every element
of $\CS(j+1,d,p)$ is the $p^{\text{th}}$ power of an element in
$\CS(j,d,p)$ (see \cite[Lemma 5.1]{DDMS99}).  It follows that
\[
\CS(j+1,d,p)=[\CS(j,d,p),\CS(j,d,p)]\CS(j,d,p)^{p}.
\]
If $t > j$, this equality implies that $\CS(j+1,d,p)/\CS(t,d,p)$ is
the Frattini subgroup of $\CS(j,d,p)/\CS(t,d,p)$.  In particular,
if we pull back $e$ ``basis'' generators of $\CS(j,d,p)/\CS(j+1,d,p)$
to $\CS(j,d,p)$ then the images of these elements
generate $\CS(j,d,p)/\CS(t,d,p)$.

Finally, if $p$ is any prime, a direct calculation shows
\[
[\CS(i,d,p),\CS(j,d,p)] \subseteq \CS(i+j,d,p)
\]
for all positive $i,j$.

We have already introduced the simplified notation $\Gamma
=\CS(u,d,p)$.  In this spirit, we set
\[
\Gamma _{i}=\CS(u+i-1,d,p)
\]
for all positive integers $i$.  The next lemma summarizes our
discussion.

\begin{Lem} \label{Lcsp}
\begin{enumerate} [\normalfont(a)]
\item \label{Lcsp1}
$\Gamma_i$ is characteristic in $\Gamma_1 = \Gamma$
for all $i$ and $\bigcap_i \Gamma_i = 1$.
\item \label{Lcsp2}
$[\Gamma_i,\Gamma_j] \subseteq \Gamma_{i+j}$ for all $i,j$.  Even
better, if $p=2$ then $[\Gamma_i,\Gamma_j] \subseteq \Gamma_{i+j+1}$.
\item \label{Lcsp3}
$\Gamma_i/\Gamma_{i+1}$ is an elementary abelian $p$-group of rank
$e=d^2$ and the $p$-power map induces an isomorphism from
$\Gamma_i/\Gamma_{i+1}$ to $\Gamma_{i+1}/\Gamma_{i+2}$ for all $i$.
\item \label{Lcsp4}
For all $i$, each element of $\Gamma _{i+1}$ is the $p^{\text{th}}$
power of an element in $\Gamma _{i}$.
\item \label{Lcsp5}
There are elements $y_{1}, \ldots , y_{e}\in \Gamma$ such that
$y_{1}\Gamma _{t}, \ldots , y_{e}\Gamma _{t}$ generate $\Gamma
/\Gamma _{t}$ for all $t$.
\end{enumerate}
\end{Lem}

\section{The Associated Graded Ring}

There is nothing mathematically new in this section.  The intended
contribution is clarity; without recasting the original argument, we
would not have made progress on the Atiyah conjecture.

Recall the notion of an $N_p$-sequence for a group $G$,
as described in \cite[\S 3.3]{Passman77}.  It is a sequence of normal
subgroups
\[
G = H_1 \supseteq H_2 \supseteq H_3 \supseteq \cdots
\]
such that $[H_i,H_j] \subseteq H_{i+j}$ and $H_i^p \subseteq H_{ip}$
for all $i,j$.  Our first step is to renumber the sequence of
congruence subgroups in Lemma \ref{Lcsp} so that it
becomes an $N_p$-sequence for $\Gamma$.  Indeed,
simply set $\Gamma_i^* =
\Gamma_{\lceil \log_p i\rceil + 1}$, where $\lceil x
\rceil$ indicates the least integer greater than or equal $x$.
Then $\Gamma_1^* = \Gamma_1$ (which we also write $\Gamma ^{*}$)
, $\Gamma_2^* = \Gamma_2$,
\dots , $\Gamma_p^* = \Gamma_2$, $\Gamma_{p+1}^* = \Gamma_3$, \dots .
The new sequence has many successive repeats; the only jumps can
occur between $\Gamma_{p^i}^*$ and $\Gamma_{p^i + 1}^{*}$.

\begin{Lem} \label{LNp}
Let $\Gamma_i^*$ be defined as above.  Then $\{\Gamma_i^*\}$ is an
$N_p$-sequence for $\Gamma$.
\end{Lem}
\begin{proof}
Since $(\lceil \log_p i\rceil + 1) + (\lceil \log_p j\rceil + 1)
\ge \lceil \log_p (i+j) \rceil + 1$, it follows from Lemma
\ref{Lcsp}\eqref{Lcsp2} that $[\Gamma_i^*,\Gamma_j^*] \subseteq
\Gamma_{i+j}^*$.  Furthermore, Lemma \ref{Lcsp}\eqref{Lcsp3} and
$\lceil \log_p ip\rceil + 1 = \lceil \log_p i\rceil + 2$ imply
$(\Gamma_i^*)^p \subseteq \Gamma_{ip}^*$.
\end{proof}

Let $t$ be a positive integer, and set $\Gamma ^{\circledast} =
\Gamma /\Gamma_t$, a finite group of order $p^{(t-1)e}$.  The
homomorphic images
$\Gamma ^{\circledast}_i = \Gamma_i^*/\Gamma_t$ for $t \ge \lceil
\log _{p} i \rceil$
constitute an $N_p$-sequence
for $\Gamma ^{\circledast}$.  Let $x_1,\dots,x_e \in
\Gamma ^{\circledast}$ denote the generators $y_{1}\Gamma _{t},
\ldots , y_{e}\Gamma _{t}$ from Lemma \ref{Lcsp}\eqref{Lcsp5}.
(These elements make up a ``basis" for
$\Gamma ^{\circledast}/\Gamma ^{\circledast}_2$.)
We apply Passman's exposition of \cite[\S
3.3]{Passman77} to the finite group $\Gamma ^{\circledast}$.  For
$g \in \Gamma ^{\circledast}$, he defines the height $\nu(g)$ to be
the greatest positive
integer $n$ such that $g \in \Gamma ^{\circledast}_n$.  By agreement,
$\nu(1) = \infty$.

Once and for all, \emph{assume that $k$ is a field of
characteristic $p$}.

Next, Passman defines $E_i
\subseteq k[\Gamma ^{\circledast}]$ to be the $k$-linear span of
all products of the form
\[
(g_1-1)(g_2-1) \dots (g_h - 1)
\]
with $\nu(g_1) + \nu(g_2) + \dots + \nu(g_h) \ge i$.
It is clear that $E_i$ is an ideal of $k[\Gamma ^{\circledast}]$ for
all $i$ and we have a descending filtration
\[
k[\Gamma ^{\circledast}] = E_0 \supseteq E_1 \supseteq E_2
\supseteq \cdots .
\]
At this point,
Passman refines the series $\{ \Gamma ^{\circledast}_{i}\}$, obtaining
a composition series. (One caveat: the new series
does not have repeated
terms even though the original series does.) He constructs a sequence
of elements in $\Gamma ^{\circledast}$ by choosing, in order, a
preimage for each generator of a composition factor.  Because of
Lemma \ref{Lcsp}\eqref{Lcsp3}, we may choose the sequence as
\[
x_1,\dots,x_e,x_1^p,\dots,x_e^p,\dots,x_1^{p^{t-2}},\dots,
x_e^{p^{t-2}}.
\]
Notice that we have $\nu(x_i^{p^s})=p^s$.  According to \cite[Lemma
3.3.5]{Passman77}, the products
\[
(x_1 - 1)^{b(1,1)} (x_2-1)^{b(2,1)} \dots
(x_{e-1}^{p^{t-2}} - 1)^{b(e-1,t-1)}
(x_e^{p^{t-2}} - 1)^{b(e,t-1)}
\]
with
\[
0 \le b(m,n) < p \text{ for all } m,n \text{ and } \sum_{i=1}^e
\sum_{j=0}^{t-2}
p^{j}b(i,j+1) \ge s
\]
constitute a basis for $E_s$.

We show that the associated graded ring
$\gr(k[\Gamma ^{\circledast}])$ is commutative.  Since
$(x_i-1)^{p^j} = x_i^{p^j} - 1$, it is
sufficient to show that $x_i-1$ and $x_j-1$ commute in $\gr(k[\Gamma
^{\circledast}])$,
or equivalently, $(x_i-1)(x_j-1) - (x_j-1)(x_i-1) \in E_3$.
Now $(x_i-1)(x_j-1) - (x_j-1)(x_i-1) = (1-[x_i,x_j])x_ix_j$, so we
need to prove $1-[x_i,x_j] \in E_3$, i.e.,
$[\Gamma , \Gamma] \subseteq \Gamma ^{*}_{3}$.  This can be read off
of Lemma \ref{Lcsp}\eqref{Lcsp2}, with a bit of care paid to the
prime $2$.  We have established that the graded ring is commutative.

The upshot of the argument so far is that $E_s$ has a better basis,
namely all products
\[
(x_1-1)^{a(1)}(x_2-1)^{a(2)} \dots (x_e-1)^{a(e)}
\]
subject to $a(j) < p^{t-1}$ for all $j$ and $\sum_{j=1}^e a(j)
\ge s$.  Obviously $E_s \subseteq \omega(\Gamma ^{\circledast})^s$.
On the other hand
$\Gamma ^{\circledast} = \langle x_1, \dots, x_e\rangle$, whence
$\omega(\Gamma ^{\circledast})
= E_1$.  Using the new basis and the commutativity of $\gr(k\Gamma
^{\circledast})$, we deduce that $E_s =
\omega(\Gamma^{\circledast})^s$ for all $s \in \mathbb{Z}_{\ge 0}$.

A quick glance at the exponent restrictions in the previous paragraph
establishes the following consequence upon which Lazard's strategy
rests.
\begin{Lem} \label{Lpolyring}
Let $\Gamma ^{\circledast}=\Gamma /\Gamma _{t}$ as above.  Then the
associated graded ring $\gr(k[\Gamma ^{\circledast}]) =
\bigoplus_{i\ge 0} \omega(\Gamma ^{\circledast})^i/\omega(\Gamma
^{\circledast})^{i+1}$ is
isomorphic to a polynomial ring in $e$ indeterminates modulo degree
$p^{t-1}$ or greater.
\end{Lem}

On a parenthetical note, recall that the $n^{\text{th}}$ dimension
subgroup for a group $G$ at the prime $p$ is
\[
D_{n}(G)= \{ g \in G \mid g-1 \in \omega (G)^{n}\}.
\]
(Here the augmentation ideal is taken inside $k[G]$ for a field $k$ of
characteristic $p$.) Observe that by definition, if $y \in \Gamma
^{\circledast}_{n}$ then $y-1 \in E_{n}$.  Since we have just proved
that $E_{n}=\omega (\Gamma ^{\circledast})^{n}$, we see that $\Gamma
^{\circledast}_{n} \subseteq D_{n}(\Gamma ^{\circledast})$.  But
$D_{n}(\Gamma ^{\circledast})$ is the smallest $N_{p}$-sequence of
$\Gamma ^{\circledast}$ \cite[Theorem 3.3.7]{Passman77}.  Therefore
$\Gamma ^{\circledast}_{n} = D_{n}(\Gamma ^{\circledast})$.  It is not
difficult to make the jump and conclude that $\Gamma ^{*}_{n}$ is the
$n^{\text{th}}$ dimension subgroup for $\Gamma ^{*}$.

\begin{Lem} \label{Lequiv}
$\omega(\Gamma)^{(p^{t-1}e-e+1)} \subseteq \omega(\Gamma_t)
k[\Gamma]\subseteq \omega(\Gamma)^{p^{t -1}}$ for all $t$.
\end{Lem}
\begin{proof}
We work inside $\Gamma ^{\circledast}=\Gamma /\Gamma _{t}$ and first
show that $\omega(\Gamma ^{\circledast})^{(p^{t-1}e-e+1)}$ is the
zero ideal.  Recall that this ideal has a basis
\[
(x_1-1)^{a(1)}(x_2-1)^{a(2)} \dots (x_e-1)^{a(e)}
\]
subject to $a(j) \le p^{t-1}-1$ and $\sum_{j=1}^e a(j)
\ge p^{t-1}e-e+1$.  Since there are no such $a(j)$, the basis is
empty, i.e., the requisite ideal is zero.

As to the second inclusion, Lemma \ref{Lcsp}\eqref{Lcsp4} implies
that each element of $\Gamma _{t}$ has the form $y^{p^{t-1}}$ for
some $y \in \Gamma$.  Hence $\omega (\Gamma _{t})$ is spanned by all
$y^{p^{t-1}}-1 = (y-1)^{p^{t-1}}$ as $y$ ranges over $\Gamma$.
\end{proof}

 From this we deduce Theorem \ref{Tcharp}.

\begin{Cor} \label{Cpolyring}
The associated graded ring $\gr(k[\Gamma]) = \bigoplus_{i=0}^{\infty}
\omega(\Gamma)^i/\omega(\Gamma)^{i+1}$
is isomorphic to the polynomial ring in $e$ variables.  Furthermore
$\bigcap \omega (\Gamma)^i = 0$ and consequently
$k[\Gamma]$ is a domain.
\end{Cor}
\begin{proof}
Initially fix $t$ and consider $\Gamma ^{\circledast}= \Gamma /\Gamma
_{t}$.  According to Lemma \ref{Lequiv}, $\omega (\Gamma
_{t})\subseteq \omega (\Gamma )^{p^{t-1}}$.  It follows that the
kernel of the map from $\omega (\Gamma )^{i}$ to $\omega (\Gamma
^{\circledast})^{i}$ coincides with the kernel of the natural map
from $k[\Gamma]$ to $k[\Gamma ^{\circledast}]$ whenever $i \le
p^{t-1}$.  Therefore $\gr(k[\Gamma])$ and
$\gr(k[\Gamma^{\circledast}])$ are isomorphic modulo degree $p^{t-1}$.
Apply Lemma \ref{Lpolyring} and let $t$ tend to infinity to obtain
the isomorphism of graded rings.

To see that $k[\Gamma]$ is a domain, we now need only check that
\[
\bigcap_i \omega (\Gamma)^i = 0.
\]
Since $\bigcap_j \Gamma_j = 1$, every nonzero element $a$ of
$k[\Gamma]$ survives in $k[\Gamma /\Gamma _{t}]$ for large enough
$t$.  That is, $a\not\in \omega (\Gamma _{t})k[\Gamma]$ for some $t$.
By Lemma \ref{Lequiv}, $a\not\in \bigcap_i \omega (\Gamma )^{i}$
\end{proof}

\section{Completed Group Rings} \label{Scomp}

The main result in this section, Theorem
\ref{Ttheorem1}\eqref{Ttheorem1b}, appears as \cite[Theorem
8.7.10]{Wilson98} and a special case of \cite[Theorem
C]{BrownArdakov05} when the field $k$ is finite.  Recall that for us,
$k$ is any field of characteristic $p$.

Let $R$ be a ring with a descending sequence $I_{0},I_{1}, I_{2},
\ldots $ of ideals.  The \emph{completion} of $R$ is
\[
\overline{R}=\varprojlim R/I_{j}.
\]
In this section, we assume that each $R/I_{j}$ has the discrete
topology.  However, $\overline{R}$ has
the usual inverse limit topology
that makes it into a Hausdorff topological ring. (This is an example
of a linear topology on $R$.) There is a more prosaic description of
$\overline{R}$ under this set-up that the reader may find useful.
Define a sequence $r_{1}, r_{2}, r_{3}, \ldots $ in $R$ to be Cauchy
provided that for every $s\ge 0$ there exists $N\ge 0$ such that
$r_{i}-r_{j}\in I_{s}$ for all $i,j\ge N$.  Then $\overline{R}$ is the
ring of all Cauchy sequences modulo the usual equivalence relation
used in advanced calculus.  For a neighborhood base at zero, we may
take all closures $\overline{I _{j}}$.

 From this description, it is clear that if $J_{0}, J_{1}, J_{2},
\ldots $ is a second descending sequence for $R$ and if for all $m,n$
there exist $m' , n'$ with
\[
I_{m}\supseteq J_{m'} \text{ and } \ J_{n}\supseteq I_{n'}
\]
then the two completed rings are isomorphic via a homeomorphism.  For
example, Lemma \ref{Lequiv} tells us that the completions of
$k[\Gamma]$ given by the descending sequences
\[
k[\Gamma] \supseteq \omega (\Gamma _{1})k[\Gamma] \supseteq
\omega (\Gamma _{2})k[\Gamma] \supseteq \cdots \text{ and }
\]
\[
k[\Gamma] \supseteq \omega (\Gamma ) \supseteq \omega (\Gamma )^{2}
\supseteq \omega (\Gamma )^{3} \supseteq \cdots
\]
are the same.  We denote this common completion by $k[[\Gamma]]$.

We will need several easy exercises whose proofs are left to the
reader.
\begin{Lem} \label{Lexer}
Assume $R=I_{0} \supseteq I_{1} \supseteq I_{2} \supseteq \cdots$ is
a descending sequence of ideals for $R$.
\begin{enumerate} [\normalfont(a)]
\item
The canonical homomorphism
$R \rightarrow \overline{R}$ induces an isomorphism
\[
R/I_{t} \rightarrow \overline{R}/\overline{I_{t}}
\]
for each $t$.
\item
$\bigcap _{s}\overline{I_{s}}=0$.
\end{enumerate}
\end{Lem}

The following lemma is undoubtedly well known to experts in
topological rings but has not been exploited in our context.

\begin{Lem} \label{Lfg}
Assume that $F$ is a commutative ring.
Let $\overline{R}$ be the completion of
the $F$-algebra $R$ with respect to a descending filtration of ideals
\[
R=I_{0} \supseteq I_{1} \supseteq I_{2} \supseteq \cdots
\]
in which $R/I_s$ is an artinian $F$-algebra for all $s$.  Then
finitely generated left (right) ideals of $\overline{R}$ are closed.
\end{Lem}
\begin{Rem*}
In the case $F$ is a field, the hypothesis that $R/I_s$ is an
artinian $F$-algebra is the same as $I_s$ having finite codimension
in $R$.
\end{Rem*}

\begin{proof}[Proof of Lemma \ref{Lfg}]
We take advantage of the theory of linearly compact modules over a
Hausdorff linearly topological ring. (A detailed discussion of this
topic can be found in \cite[Chapter VII]{Warner93}.)

We consider $R/I_{n}$ as a left (or right) $\overline{R}$-module by
Lemma \ref{Lexer}.  Since $R/I_{n}$ is artinian over $F$, it
is linearly compact \cite[28.14]{Warner93}.  According to the argument
in \cite[28.15]{Warner93} (essentially Tychonoff's theorem), the
inverse limit of linearly compact modules is linearly compact.
Hence $\overline{R}$ is linearly compact as a module over itself.

If $a \in \overline{R}$ then $\overline{R}a$ is linearly compact
because it is the continuous homomorphic image of a linearly compact
$\overline{R}$-module \cite[28.3]{Warner93}.  It follows that
$\overline{R}a$ is closed in $\overline{R}$ \cite[28.6(2)]{Warner93}.
Finally, the sum of two closed submodules of a linearly compact
module is closed \cite[28.6(3)]{Warner93}.  We conclude that any
finitely generated left ideal of $\overline{R}$ is closed.
\end{proof}

The streamlined notation $\overline{\omega}(\Gamma )$ will be used for
the closure of the augmentation ideal in $k[[\Gamma]]$.

\begin{Thm} \label{Ttheorem1}
\begin{enumerate}[\normalfont(a)]
\item \label{Ttheorem1a}
The associated graded ring $\gr(k[[\Gamma]])=\bigoplus_{i=0}^{\infty}
\overline{\omega}(\Gamma)^i/\overline{\omega}(\Gamma)^{i+1}$
is isomorphic to the polynomial ring in $e$ variables.

\item \label{Ttheorem1b}
The completed group ring $k[[\Gamma]]$ is a noetherian domain.
\end{enumerate}
\end{Thm}
\begin{proof}
Recall that $\bigcap_{j}\omega (\Gamma )^{j} =0$ by Corollary
\ref{Cpolyring}.  It then follows from Lemma \ref{Lexer} that
the powers of $\overline{\omega}(\Gamma )$ make up a
descending filtration for $k[[\Gamma]]$ with
\[
\gr k[[\Gamma]] \simeq \gr k[\Gamma].
\]
Part \eqref{Ttheorem1a} now follows from Corollary \ref{Cpolyring}.
Furthermore, since $\gr k[[\Gamma]]$ is
a domain, we conclude that the completion $k[[\Gamma]]$ is, too.

To see that $k[[\Gamma]]$ is (left) noetherian, consider a nonzero
left ideal $L$.  The collection of its bottom symbols in
$\gr k[[\Gamma]]$ constitute a homogeneous left ideal, $\gr L$,
inside the graded algebra.  But $\gr k[[\Gamma]]$ is
isomorphic to a polynomial
ring in $e$ indeterminates and thus is noetherian.  Choose finitely
many homogeneous generators for $\gr L$ and lift them to elements
$b_{1}, b_{2}, \ldots , b_{m}$ in $L$.  Suppose $y\in L$ is nonzero
and has degree $s$.  We peel off the bottom of $y$: there exist
$r_{1}, \ldots , r_{m}\in k[[\Gamma]]$ such that the bottom symbol
of
\[
y-r_{1}b_{1} - \cdots - r_{m}b_{m}
\]
has degree greater than $s$.  Continuing inductively in this fashion,
we see that $y$ lies in the closure of the left ideal generated by
$b_{1}, \ldots , b_{m}$.  However, this ideal is already closed
(apply Lemma \ref{Lfg} with $F=k$).
Therefore $L$ itself is finitely generated.
\end{proof}

\section{A Key Estimate} \label{Skey}

Our goal for the rest of the paper is to settle the Atiyah problem for
virtual congruence subgroups.  The strategy is to establish the
algebraic ingredient in L\"uck's approach: an approximate nullity of
matrices over the group algebra always takes on integer values.  We
accomplish this in a series of steps.  In this section, we address
the ``one-dimensional case'' for $k[[\Gamma]]$.  The next step is to
extend this estimate ``up by finite index".  At the end, we apply
gaussian elimination to reduce the matrix computation to the scalar
case.

We continue using the notation $\Gamma _{i} = \CS(u+i-1,d,p)$ with
$\Gamma = \Gamma _{1}$.  As usual, $k$ is a field of
characteristic $p$.

\begin{Lem} \label{Lprekey}
Fix a positive integer $m$ and let $\pi \colon k[[\Gamma]] \to
k[\Gamma]/\omega (\Gamma )^{m}$ be the canonical homomorphism.
Assume $0\neq a \in k[[\Gamma]]$ has degree $s$ where
$m>s$.  Denote by $L \colon \pi (k[[\Gamma]]) \to \pi
(k[[\Gamma]])$ the map ``multiply on the left by $\pi (a)$".  Then
\[
\ker L = \pi (\overline{\omega}(\Gamma )^{m-s}).
\]
\end{Lem}
\begin{proof}
Since $a\in \pi (\overline{\omega}(\Gamma )^{s})$, we have
\[
\pi (\overline{\omega}(\Gamma )^{m-s}) \subseteq \ker L.
\]
On the other hand, suppose $0\neq b \in k[[\Gamma]]$ with $\pi (b)
\in \ker L$.  Let $t$ be the degree of $b$.  According to Theorem
\ref{Ttheorem1}\eqref{Ttheorem1a},
$ab \not\in \overline{\omega}^{s+t+1}$.  But $\pi
(ab)=0$, so $ab \in \overline{\omega}(\Gamma )^{m}$.
Thus $s+t+1 > m$.  In other words, $t \ge m-s$.  Therefore
\[
\pi (b) \in \pi (\overline{\omega}(\Gamma )^{t})\subseteq \pi
(\overline{\omega}(\Gamma )^{m-s}).
\]
\end{proof}
Recall that there is a ring homomorphism $k[[\Gamma]] \to
k[\Gamma /\Gamma _{n}]$ for all positive integers $n$ by Lemma
\ref{Lexer}.  Use the notation $[\Gamma \colon \Gamma _{n}]=
|\Gamma /\Gamma _{n}|$.

\begin{Lem} \label{Lkey}
Fix $0 \ne a \in k[[\Gamma]]$ and for each $n\ge 1$ let
$\alpha_n \colon k[\Gamma/\Gamma_n]\to k[\Gamma/\Gamma_n]$
denote the map induced from left multiplication by $a$ in
$k[[\Gamma]]$.  Then
\[
\lim_{n\to \infty} (\dim_k \ker\alpha_n)/[\Gamma \colon \Gamma_n] = 0.
\]
\end{Lem}
\begin{proof}
Initially fix $n$ and set $m=p^{n-1}e-e+1$.  By
Lemma \ref{Lequiv}, $\omega (\Gamma )^{m}
\subseteq \omega (\Gamma _{n})k[\Gamma]$.  Hence there
is a natural right $k[[\Gamma]]$-module surjection $k[\Gamma
]/\omega (\Gamma )^{m}\to k[\Gamma]/\omega (\Gamma _{n})k[\Gamma]$.
Let
\[
L_{n} \colon k[\Gamma]/\omega (\Gamma )^{m}\to k[\Gamma]/\omega
(\Gamma )^{m}
\]
be the map induced from left multiplication by $a$.  Since this
multiplication map is a right $k[[\Gamma]]$-module map,
\[
\dim_k \ker L_n = \dim_k \coker L_n \ge \dim_k \coker\alpha_n =
\dim_k\ker\alpha_n.
\]
Thus it will be sufficient to show that $\lim_{n \to \infty} (\dim_k
\ker L_n)/[\Gamma \colon \Gamma_n] = 0$.

Let $s$ be the degree of $a$ in $k[[\Gamma]]$.
Lemma \ref{Lprekey} established that
\[
\dim_k \ker L_n = \dim_k (\omega
(\Gamma )^{m-s}/\omega (\Gamma )^{m}).
\]
This number is the same as the number of monomials in the polynomial
ring with $e$ variables each of which has total degree in the interval
$[m-s , m)$.  It is well known that the number of monomials of degree
less than $r$ is $\binom{e+r-1}{e}$.  Hence
\[
\dim_k \ker L_n = \binom{e+m-1}{e}-\binom{e+m-s-1}{e}.
\]
As a polynomial in $m$, each binomial coefficient has leading term
$m^{e}$.  Their difference is a polynomial in $m$ with degree at most
$e-1$.  Thus there is a constant $C$ (independent of $m$) such that
\[
\dim_k \ker L_n \le C(m+e-1)^{e-1}.
\]
Now $[\Gamma \colon \Gamma _{n}] = p^{e(n-1)}$.  Therefore
\[
\frac{\dim_k \ker L_n }{[\Gamma \colon \Gamma _{n}]} \le
C(p^{n-1}e)^{e-1}p^{-e(n-1)}= \frac{Ce^{e-1}}{p^{n-1}}.
\]
Take the limit as $n$ tends to $\infty$.
\end{proof}

\section{Virtual Congruence Subgroups} \label{Svirt}

Henceforth, we will reserve the symbol $G$ for a torsion-free
group that contains a normal subgroup of finite index isomorphic to
$\Gamma$.  By abuse of notation, we identify this subgroup with
$\Gamma$.

\begin{Thm} \label{TOre}
$k[[G]]$ is a right and left noetherian domain.
\end{Thm}
\begin{proof}
Clearly $G$ is a profinite group, and since it is torsion-free, it
is easy to see from \cite[Proposition 2.3.3]{Wilson98} that $G$ is
a pro-$p$ group.  A version of the construction discussed in
Section \ref{Scomp} makes it legitimate to write $k[[G]]$.

We know by Theorem \ref{Ttheorem1}\eqref{Ttheorem1b}
that $k[[\Gamma]]$ is a noetherian domain.
Let $S$ denote the nonzero elements of $k[[\Gamma]]$.  Then $S$ is
a right Ore set in $k[[\Gamma]]$ that does not contain zero divisors.
This means the ring of fractions $k[[\Gamma]]S^{-1}$
is a division ring.  This ring consists
of elements of the form $as^{-1}$ with $a \in
k[[\Gamma]]$ and $s \in S$; however because $k[[\Gamma]]$ is left
noetherian, it also consists of elements of the form $s^{-1}a$
with $a \in k[[\Gamma]]$ and $s \in S$.
By \cite[Proposition 7.6.3]{Wilson98},
$k[[G]]$ is a free right $k[[\Gamma]]$-module with basis any (left)
transversal for $\Gamma$ in $G$.  Since $G/\Gamma$ is
finite, we see that $k[[G]]$ is a free right
$k[[\Gamma]]$-module in the usual module theoretic sense.
It now follows that $k[[G]]$ is noetherian and that
$S$ is also a right Ore set in
$k[[G]]$, so we can form the
ring $k[[G]]S^{-1}$.  This will be an Artinian ring and
a right $k[[\Gamma]]S^{-1}$-module containing $k[[G]]$.

We know by \cite[Proposition
8.5.1]{Wilson98} that $\Gamma$ is a pro-$p$ group of finite rank.
Therefore by \cite[Proposition 8.1.1]{Wilson98}, $G$ is a profinite
group of finite rank; we deduce from \cite[Theorem
11.6.9]{Wilson98} that $G$ has finite cohomological dimension.
It now follows from the results of \cite{Brumer66} that $k[[G]]$ has
finite global dimension.  In fact, on \cite[p.~443]{Brumer66} of
Brumer's paper, one finds
that if $G$ is a pro-$p$ group of finite cohomological
dimension and $\Omega$ is a complete regular local ring in which $p$
is not a unit then $\Omega[[G]]$ is a complete noncommutative
local ring of finite global dimension.  We apply this with $k=
\Omega$.

Suppose $k[[G]]S^{-1}$ is not a division ring.  We borrow the
strategy of using Walker's Theorem \cite{Walker72}
as first introduced by Ken Brown \cite{Brown76}.
There must be an ideal $I$ of $k[[G]]S^{-1}$
such that $0 < \dim_{k[[\Gamma]]S^{-1}} I < [G:\Gamma]$.
Set $J = k[[G]] \cap I$, so $I$ is a right ideal of $k[[G]]$
and $JS^{-1} = I$.
Since $k[[G]]$ is a right noetherian ring with finite cohomological
dimension, there is a resolution
\[
0 \longrightarrow P_n \longrightarrow P_{n-1} \longrightarrow \dots
\longrightarrow P_0 \longrightarrow J \longrightarrow 0
\]
where the $P_i$ are finitely generated projective $k[[G]]$-modules.
By \cite[Corollary 7.5.4]{Wilson98}, each $P_i$ is a free
$k[[G]]$-module, and since $P_i$ is finitely generated, this simply
means that $P_i$ is free in the usual module theoretic sense.
Applying $S^{-1}$ to the resolution of $J$, we obtain a resolution of
$JS^{-1}$ with free $k[[G]]S^{-1}$-modules.  We deduce that
$\dim_{k[[\Gamma]]S^{-1}} JS^{-1}$ is a multiple of $[G:\Gamma]$.
This is a contradiction.
\end{proof}

Suppose that $R$ is a complete discrete valuation ring with maximal
ideal $qR$.  Assume, in addition, that $R/qR$ is a field of
characteristic $p$.  We regard $R$ as a Hausdorff topological
ring as in Section \ref{Scomp} by assuming that each $R/q^jR$
has the
discrete topology and then giving $R$ the inverse limit topology.
Now define $R[[G]] = \varprojlim R[G/N]$, where the inverse
limit is taken over the open normal subgroups $N$ of $G$.
(Here each $R[G/N]$ is given the product
topology of a finitely generated free $R$-module.) By construction,
$\bigcap_n q^n R[[G]] = 0$.  It is also easy to see that nonzero
members of $R$ cannot become zero divisors in $R[[G]]$.

\begin{Cor} \label{COre}
Let $R$ be a complete discrete valuation ring
with maximal ideal $qR$ such that $R/qR$ is a field of
characteristic $p$.  Then $R[[G]]$ is a domain.
\end{Cor}
\begin{proof}
Set $k=R/qR$.  There is a continuous surjective ring homomorphism
\[
\theta \colon R[[G]] \to k[[G]]
\]
with kernel $qR[[G]]$.  Suppose $0 \neq a,c \in R[[G]]$.  Write
\[
a=q^{m }\widetilde{a} \quad\text{and}\quad c=q^{n}\widetilde{c}
\]
with $\widetilde{a}, \widetilde{c}\not\in qR{{G}}$.  Then
\[
\theta (\widetilde{a}\,\widetilde{c})=\theta (\widetilde{a})\theta
(\widetilde{c}) \neq 0
\]
by Theorem \ref{TOre}.  We conclude that $\widetilde{a}\widetilde{c}
\neq 0$, so $ac \neq 0$.
\end{proof}

Since $\Gamma _{n}$ is characteristic in $\Gamma$, each of these
subgroups is normal in $G$.  It follows from
\cite[Proposition 7.1.2(c)]{Wilson98} that there is a homomorphism
$k[[G]] \rightarrow k[[G/\Gamma _{n}]]$ for all $n$.  In particular,
left multiplication by an element of $k[[G]]$ induces a right
module endomorphism of $k[G/\Gamma _{n}]$ for each $n$.

\begin{Lem} \label{Lextension}
Suppose $a \in k[[G]]$ is nonzero and $\alpha_n \colon
k[G/\Gamma_n] \to k[G/\Gamma_n]$ denotes the map induced from left
multiplication by $a$ in $k[[G]]$.  Then
\[
\lim_{n\to \infty} (\dim_k \ker\alpha_n)/[G \colon \Gamma_n] = 0.
\]
\end{Lem}
\begin{proof}
By Theorem \ref{TOre}, $k[[G]]$ is a domain and $k[[\Gamma]]$
is a noetherian Ore domain inside.  As we have already observed in
the proof of Theorem \ref{TOre},
if $S = k[[\Gamma]]\setminus 0$ then we may form the division
ring of fractions $S^{-1}k[[G]]$.  This means that we can find
$b \in k[[G]]$ such that $0 \ne ba \in k[[\Gamma]]$.
For each non-negative integer $n$, let $\sigma_n \colon
k[G/\Gamma_n] \to k[G/\Gamma_n]$ denote the right
$k[G/\Gamma_n]$-map induced from left multiplication by $ba$ on
$k[[G]]$.  Similarly, let $\sigma_n |\, \colon k[\Gamma /\Gamma_n]
\to k[\Gamma /\Gamma_n]$ be the corresponding map for the restriction
of multiplication by $ba$ to the subalgebra $k[[\Gamma]]$.
Clearly
\[
\dim_{k} \ker\sigma _{n}=[G\colon \Gamma]\dim_{k} \ker(\sigma
_{n}|\, ).
\]
Since $[G\colon \Gamma _{n}]=[G\colon \Gamma][\Gamma \colon \Gamma
_{n}]$, Lemma \ref{Lkey} implies
\[
\lim_{n \to \infty} (\dim_k
\ker\sigma_n)/[G\colon\Gamma_n] = 0.
\]
But $\dim_{k}\ker\alpha_n \le \dim_k \ker\sigma_n$.
We deduce that the desired limit is zero.
\end{proof}

For the rest of this paper, we will assume that $k$
is a \emph{finite} field
of characteristic $p$.  This choice is mainly made for convenience; we
would like to quote Wilson \cite{Wilson98}, who requires that his
completed algebras be profinite.  Mathematically, this is primarily
an issue of applying compactness rather than linear compactness.
We can now state
\begin{Prop} \label{POre}
Let $R$ be a discrete valuation ring with maximal ideal $qR$ such
that $R/qR$ is a finite field of characteristic $p$.  Then $R[[G]]$
is a noetherian domain.
\end{Prop}
\begin{proof}
According to
\cite[Propositions 8.1.1, 8.5.1 and Theorem 8.7.8]{Wilson98},
$R[[G]]$ is noetherian.  The result now follows from Corollary
\ref{COre}.
\end{proof}

\begin{Rem*}
Actually Proposition \ref{POre} remains true if $R/qR$ is an arbitrary
field of characteristic $p$.  We sketch the argument for this.  We
already know that $R[[G]]$ is a domain by Corollary \ref{COre}, so it
remains to prove that $R[[G]]$ is noetherian.  To establish this,
it will be sufficient to descend to a subgroup of finite index in $G$.
We may now assume that $G$ is ``extra powerful" \cite[p.~148 and
Proposition 8.5.2]{Wilson98}.  Let $\Delta$ denote the unique maximal
ideal of $R[[G]]$.  Then $\gr(R[[G]]) := \bigoplus_{i=0}^{\infty}
\Delta^i/\Delta^{i+1}$ is a commutative noetherian ring
by the proofs of \cite[Theorems 8.7.6 and 8.7.7]{Wilson98}.  We can
now use Lemma \ref{Lfg} and the proof of Theorem
\ref{Ttheorem1}\eqref{Ttheorem1b} to establish that $R[[G]]$ is
noetherian.
\end{Rem*}

Suppose that $R$ is a ring as described in the previous corollary.
Since $R$ is a principal ideal domain, any submodule $M$ of a
finitely generated free $R$-module is again free. (During the
discussion that follows, the free module will be the ring of all
matrices of some size over $R[H]$ for a finite group $H$.) In
particular, $M$ has a well-defined rank, which we perversely write
as $\dim _{R}M$ to conform with earlier notation.  Note that if $K$
is the field of fractions for $R$ then
\[
\dim _{K}K\otimes _{R}M=\dim _{R}M.
\]

\begin{Cor} \label{Cextension}
Let $R$ be a complete discrete valuation ring
with maximal ideal $qR$ such that $R/qR$ is a finite field of
characteristic $p$.
Suppose $a \in R[[G]]$ is nonzero and $\alpha_n \colon
R[G/\Gamma_n] \to R[G/\Gamma_n]$ denotes the map induced from left
multiplication by $a$ in $R[[G]]$.  Then
\[
\lim_{n\to \infty} (\dim_R \ker \alpha_n)/[G:\Gamma_n] = 0.
\]
\end{Cor}
\begin{proof}
Write $a = q^fc$ where the non-negative integer $f$ is chosen so that
$c \in R[[G]]\setminus qR[[G]]$.  Let $\xi_n
\colon R[G/\Gamma_n] \to R[G/\Gamma_n]$ be the multiplication map
induced from $c$.  If $k$ denotes the residue field for $R$, we obtain
a right $k[G/\Gamma _{n}]$-endomorphism by factoring out $q$,
\[
\xi ^{*}_{n} \colon k[G/\Gamma _{n}] \to k[G/\Gamma _{n}].
\]
(It is induced from left multiplication by the \emph{nonzero} image
of $c$ in $k[[G]]$.)

We may choose an $R$-basis for $\ker \xi _{n}$ that
consists of elements in $R[G/\Gamma _{n}]\setminus qR[G/\Gamma _{n}]$.
Their images modulo $q$ remain linearly independent.  Hence
\[
\dim _{R}\ker \xi _{n} \le \dim _{k}\ker \xi ^{*}_{n}.
\]
But $\dim _{R} \ker \alpha _{n}=
\dim _{R} \ker \xi _{n}$.  The required inequality is now a
consequence of Lemma \ref{Lextension}.
\end{proof}

\section{Matrices}

\begin{Lem} \label{Lmatrix}
Let $R$ be a complete discrete valuation ring with maximal ideal
$qR$ such that $R/qR$ is a finite field of characteristic $p$, and
let $h$ be a positive integer.
For $A \in \Mat_h(R[[G]])$ and for each
positive integer $n$, let $\mathcal{A}_n \colon R[G/\Gamma_n]^h \to
R[G/\Gamma_n]^h$ denote the right $R[G/\Gamma_n]$-module map
induced from acting by $A$ on the left.  Then $\lim_{n\to \infty}
(\dim_R \ker \mathcal{A}_n)/[G: \Gamma_n]$ is an integer.
\end{Lem}
\begin{proof}
Since $R[[G]]$ is an Ore domain by Proposition \ref{POre},
it has a division ring of fractions $Q$.  By performing complete
row and column reduction, we can find invertible $B,C \in \Mat_h(Q)$
such that $BAC$ is a diagonal matrix $D$ with only zeros and
ones on the main diagonal.  Let
$\delta$ denote the number of such zeros in $D$
(so $0 \le \delta \le h$).  Using the fact that $R[[G]]$ is a left
and right Ore domain, we can find nonzero $b,c$ in $R[[G]]$ such that
$bB,Cc \in \Mat_h(R[[G]])$.  Now $(bB)A(Cc)$ is a diagonal matrix
with $bc$ appearing $h-\delta$ times on the diagonal and zero
appearing $\delta$ times down the diagonal.  Let
$\alpha_{n} \colon R[G/\Gamma_n] \to R[G/\Gamma_n]$
denote the right $R[G/\Gamma_n]$-map induced from left
multiplication by $bc$.  Similarly, let $\mathcal{B}_n,
\mathcal{C}_n \colon R[G/\Gamma_n]^h \to R[G/\Gamma_n]^h$ be the
module maps induced from the left action of $bB$ and $Cc$
respectively.  Then
\[
\dim_R \ker \mathcal{A}_n \le \dim_R \ker (\mathcal{B}_n\mathcal{A}_n
\mathcal{C}_n) = (h-\delta )\dim _{R}\ker\alpha _{n} + \delta
[G: \Gamma _{n}].
\]
But $\lim_{n \to \infty} \dim_R \ker
\alpha_{n} /[G: \Gamma_n] = 0$ by Corollary \ref{Cextension}.
We conclude that
\begin{equation} \label{Ematrix1}
\varlimsup_{n \to \infty} (\dim_R \ker \mathcal{A}_n)/[G:\Gamma_n]
\le \delta.
\end{equation}

Let $E \in \Mat_n(Q)$ denote the diagonal matrix obtained from $D$ by
replacing each one with zero and each zero on the diagonal with one.
The upshot is that $E$ has $h-\delta$ diagonal zeros and
$DE = 0$.  Since $B$ is invertible, we deduce
that $ACE = 0$.  As above, we can find a nonzero $x$ in $R[[G]]$ such
that $CEx$ is a matrix $X$ in $\Mat_h(R[[G]])$.  Then $AX = 0$ and
$C^{-1}Xx^{-1}=E$.  Choose a nonzero $c'\in R[[G]]$ so that
$c'C^{-1}\in \Mat_h(R[[G]])$.
As a consequence, $(c'C^{-1})X$ is a diagonal matrix with $c'x$
appearing $\delta$ times on the diagonal and zero appearing $h-\delta$
times down the diagonal.  Repeat the argument of the first paragraph
replacing $A$ with $X$.  We conclude that if
$\mathcal{X}_n \colon R[G/\Gamma_n]^h \to R[G/\Gamma_n]^h$ denotes the
map induced from the left action of $X$ then
\[
\varlimsup_{n \to \infty} (\dim_R \ker \mathcal{X}_n)/[G:
\Gamma_n] \le h - \delta .
\]
Therefore
$\varliminf_{n \to \infty} (\dim_R \im \mathcal{X}_n)/[G: \Gamma_n]
\ge \delta$.  We conclude from $\mathcal{A}_{n}\mathcal{X}_{n}=0$ that
\begin{equation} \label{Ematrix2}
\varliminf_{n \to \infty} (\dim_R \ker \mathcal{A}_n)/[G:
\Gamma_n] \ge \delta.
\end{equation}
Combining \eqref{Ematrix1} and \eqref{Ematrix2},
we deduce that
$\lim_{n\to \infty} (\dim_R\ker \mathcal{A}_n)/[G: \Gamma_n]
= \delta$.  Of course, $\delta$ is by definition an integer.
\end{proof}

We can finally settle the Atiyah problem for virtual congruence
subgroups.  The remaining argument essentially explains how the
integer limit lemma we have just proved fits into well established
work on this problem.

\begin{proof}[Proof of Theorem \ref{Tmain}]
An arbitrary $\mathcal{N}(G)$-module $M$ has a well-defined von
Neumann dimension $\dim_{\mathcal{N}(G)} M$ \cite[\S 6.1]{Lueck02}.
According to the proof of \cite[Lemma 10.39]{Lueck02}, it suffices to
show that for any $A \in \Mat_h(\overline{\mathbb{Q}}G)$, the induced
left action map $\mathcal{A} \colon L^2(G)^h \to L^2(G)^h$ satisfies
the integrality condition
\[
\dim_{\mathcal{N}(G)}\ker\mathcal{A} \in \mathbb {Z}.
\]
(Technically, L\"uck's proof is stated
for $\mathbb{C}$ rather than $\overline{\mathbb{Q}}$, but the proof
works just as well in this situation.) Since $A$ has only finitely
many entries, we may assume that $A \in \Mat_h(K[G])$ for some
algebraic number field $K$.  L\"uck's condition will follow if we can
show that
\[
\lim_{n\to \infty} (\dim_K \ker\mathcal{A}_{n})/[G:\Gamma_n]
\in \mathbb {Z}
\]
where $\mathcal{A}_n \colon K[G/\Gamma_n]^h \to K[G/\Gamma_n]^h$ is
the right module map induced from the left action of $A$
(see \cite[Theorem 1.6]{DLMSY03}).

Extend the $p$-adic valuation on $\mathbb{Q}$
to $K$ and complete $K$ to obtain $\hat{K}$.  Let $R$ be the
complete discrete valuation ring determined by the extended valuation.
Choose a nonzero $r\in R$ so that
$B=rA$ lies in $\Mat_{h}(R[G])$.  Define
$\mathcal{B}_{n}$ on $R[G/\Gamma_n]^h$ in the usual way from the left
action of $B$.  By the remarks preceding Corollary \ref{Cextension},
\[
\dim _{R} \ker \mathcal{B}_{n}=\dim _{\hat{K}}\ker (1\otimes
_{\hat{K}}\mathcal{A}_{n})=\dim _{K}\ker\mathcal{A}_{n}.
\]
Apply Lemma \ref{Lmatrix} to obtain $\lim_{n\to \infty}
\dim_R (\ker \mathcal{B}_n)/[G:\Gamma_n]
\in \mathbb{Z}$ and substitute.
\end{proof}

Finally we remark that the main result of this paper remains true
when $G$ is a torsion-free $p$-adic analytic (or Lie) group
\cite[Definition 8.14]{DDMS99} or, equivalently, a torsion-free
pro-$p$ group of finite rank \cite[Corollary 8.34]{DDMS99}.
Such a group has a ``uniform" open characteristic subgroup
\cite[Theorem 8.5.3]{Wilson98}, which plays the r\^ole of our
congruence subgroups; in particular one needs to replace Theorem
\ref{Ttheorem1}\eqref{Ttheorem1a} with \cite[Theorem
8.7.10]{Wilson98}.  For completeness, we state
\begin{Thm}
Let $G$ be a torsion-free pro-$p$ group of finite rank.  Then there
is a division ring $D(G)$ such that $\overline{\mathbb{Q}}[G]
\subseteq D(G) \subseteq \mathcal{U}(G)$.
\end{Thm}

\end{document}